\documentclass[12pt]{article}
\usepackage[utf8]{inputenc}
\usepackage{amsfonts}
\usepackage{scrextend}
\usepackage{mathtools}
\usepackage{float}
\usepackage{amsmath}
\usepackage[margin=1.2in]{geometry}
\usepackage{authblk}
\usepackage{amssymb}
\usepackage{epigraph}
\usepackage[english]{babel}
\usepackage[nottoc]{tocbibind}
\usepackage{tikz-cd}
\usepackage{graphicx}
\usetikzlibrary{matrix}
\usepackage{amsthm}
\usepackage[mathscr]{euscript}
 \let\mathscr\relax
\usepackage[scr]{rsfso}
\usepackage{comment}

\newtheorem{nummer}{ }[section]

\newtheorem{thm}[nummer]{\sc Theorem}
\newtheorem{prp}[nummer]{\sc Proposition}

\newtheorem{cor}[nummer]{\sc Corollary}
\newtheorem{fct}[nummer]{\sc Fact}

\newtheorem{rmk}{\sc Remark}

\newenvironment{claim}[1]{\par\noindent\textbf{Claim.}\space#1}{}

\newcounter{faelle} 
\renewcommand{\thefaelle}{\rm(\arabic{faelle})}

\newcommand\modM{\boldsymbol{M}}

\renewcommand\qed{\relax\ifmmode~\hfill$\dashv$\else\unskip\nobreak~\hfill$\dashv$\fi}

\def\epsilon{\varepsilon}

\newcommand\fin{\operatorname{fin}}
\newcommand\fix{\operatorname{fix}}

\newcommand\sym{\operatorname{sym}}

\newcommand\ZF{\text{\sf ZF}}

\newcommand\proves{\vdash}

\newcommand{\blcb} {\big{\{}}
\newcommand{\brcb} {\big{\}}}

\newcommand\subs{\subseteq}

\renewcommand{\phi}{\varphi}
\renewcommand{\theta}{\vartheta}

\newcommand{\T}[1]{\text{\sf T}_{\!#1}}
\newcommand{\RC}[1]{\operatorname{RC}_{#1}}
\newcommand{\Cm}[1]{\operatorname{C}^{-}_{#1}}
\newcommand{\RCmn}[2]{\RC{#1}\Rightarrow\RC{#2}}

\newcommand{\MOD}[1]{\text{\bf{MOD}}_{#1}}
\newcommand{\Sel}[1]{\text{\sf Sel}_{#1}}

\begin{document}

\begin{center}
{\Large\sc Implications of Ramsey Choice Principles in $\ZF$}\\[1.8ex]
{\small Lorenz Halbeisen,  
Riccardo Plati, and Saharon Shelah
\footnote{
Research partially supported by the {\it Israel 
Science Foundation\/} grant no.\;1838/19.
Paper\;1243 on author's publication list.
}}
\end{center}


\begin{quote}
{\small {\bf Abstract.} 
The Ramsey Choice principle for families of
$n$-element sets, denoted $\RC n$, states that every infinite
set $X$ has an infinite subset 
$Y\subs X$ with a choice function on $[Y]^n := \{z\subs Y : |z| = n\}$. 
We investigate for which positive integers $m$ and $n$ the implication
$\RC mn$ is provable in~$\ZF$. It will turn out that
beside the trivial implications $\RC m$, 
under the assumption that every odd integer $n>5$ 
is the sum of three primes (known as ternary Goldbach conjecture),
the only non-trivial implication which is provable in $\ZF$ is 
$\RCmn 24$.}
\end{quote}

\begin{quote}
\small{{\bf key-words\/}: permutation models, consistency results,
Ramsey choice, ternary Goldbach conjecture}\\
\small{\bf 2010 Mathematics Subject Classification\/}: {\bf 03E35}\ {03E25}
\end{quote}


\setcounter{section}{0}
\section{Introduction}

For positive integers $n$, the Ramsey Choice principle for 
families of $n$-element sets, denoted $\RC n$, is defined as follows:
For every infinite set $X$ there is an infinite subset 
$Y\subs X$ such that the set $[Y]^n := \{z\subs Y : |z| = n\}$
has a choice function. The Ramsey Choice principle was introduced 
by Montenegro~\cite{montenegro} who
showed that for $n=2,3,4$, $\RC n\Rightarrow\Cm n$. 
where $\Cm n$ is the statement that every infinite
family of $n$-element has an infinite subfamily with a choice 
function.  However, the question of whether or not 
$\RC{n}\rightarrow\mathsf{C}_{n}^{-}$ for 
$n\ge 5$ is still open (for partial answers to this question
see~\cite{LorSch,LorLeft}).

In this paper, we investigate the relation between $\RC n$ and $\RC m$
for positive integers $n$ and $m$. First, for each positive integer $m$
we construct a permutation models $\MOD m$ in which $\RC m$,
and then we show that $\RC n$ fails in $\MOD m$ for certain integers~$m$.
In particular, assuming the ternary Goldbach conjecture, which states 
that every odd integer $n>5$ is the sum of three primes, and 
by the transfer principles of Pincus~\cite{pincus}, we 
we obtain that for $m,n\ge 2$, the implication $\RCmn mn$ is 
not provable in $\ZF$ except in the case when $m=n$, or when $m=2$
and~$n=4$.

\begin{fct}\label{fct:RC24}
The implications\/ $\RCmn mm$ (for\/ $m\ge 1$) 
and\/ $\RCmn 24$ are provable in~$\ZF$.
\end{fct}

\begin{proof}
The implication $\RCmn mm$ is trivial. To see that $\RCmn 24$
is provable in $\ZF$, we assume $\RC 2$. If $X$ is an infinite
set, then by $\RC 2$ there is an infinite subset $Y\subs X$
such that $[Y]^2$ has a choice function $f_2$. Now,
for any $z\in [Y]^4$, $[z]^2$ is 
a $6$-element subset of $[Y]^2$, and by the choice function
$f_2$ we can select an element from each $2$-element subset
of~$z$. For any $z\in [Y]^4$ and each $a\in z$, 
let $\nu_z(a):=|\{x\in [z]^2:f_2(x)=a\}$, 
$m_z:=\min\blcb\nu_z(a):a\in z\brcb$, and 
$M_z:=\blcb a\in z:\nu_z(a)=m_z\brcb$.
Since $f_2$ is a choice function, we have $\sum_{a\in z}\nu_z(a)=6$,
and since $4\nmid 6$, the function $f:[Y]^4\to Y$ defined by
stipulating
$$f(z):=
\begin{cases}
a & \text{if $M_z=\{a\}$},\\
b & \text{if $z\setminus M_z=\{b\}$},\\
c & \text{if $|M_z|=2$ and $f_2(M_z)=c$},
\end{cases}
$$
is a choice function on $[Y]^4$, which shows that $\RC 4$ holds.
\end{proof}

\section{A model in which $\boldsymbol{\RC m}$ holds}

In this section we construct a permutation model $\MOD m$ in which
$\RC m$ holds. According to \cite[p.\,211\,{\sl ff.}]{cst},
the model $\MOD m$ is a {\it Shelah Model of the Second Type}.

Fix an integer $m\ge 2$ and let $\mathcal{L}_m$ 
be the signature containing the relation symbol $\Sel m$. 
Let $\T m$ be the $\mathcal{L}_m$-theory containing the 
following axiom-schema:\vspace{1em}

\begin{addmargin}[27pt]{27pt}
\textit{For all pairwise different $x_1,\dots,x_m$, 
there exists a unique index $i\in\{1,\dots,m\}$ such that, 
whenever $\{b_1,\dots,b_m\}=\{1,\dots,m\}$, 
$$\Sel m(x_{b_1},\dots,x_{b_m},x_{b})\iff b = i.$$}
\end{addmargin} \vspace{1em}

\noindent In other words, $\Sel m$ is a 
selecting function which selects
an element from each $m$-element set~$\{x_1,\ldots,x_m\}$.
In any model of the theory $\T m$, the relation 
$\Sel m$ is equivalent to a function $\Sel{}$ 
which selects a unique element from any $m$-element set.\newline

\noindent For a model $\modM$ of $\T m$ with domain $M$, 
we will simply write $M\models\T m$. 
Let $$\widetilde{C}=\{M:M\in\fin(\omega)\wedge M\models\T m\}.$$ 
Evidently $\widetilde{C}\neq\emptyset$. Partition $\widetilde{C}$ 
into maximal isomorphism classes and let $C$ be a set of 
representatives. We proceed with the construction of the set of 
atoms for our permutation model. The next theorem and its 
proof are taken from \cite{cst}, with a minor difference 
which will play an essential role in our work.

\begin{prp}\label{prp:2.1}
Let $m\in\omega\setminus \{0\}$. There exists a model 
$\mathbf{F}\models\T m$ with domain $\omega$ such that
\begin{itemize}
    \item Given a non empty $M\in C$, $\mathbf{F}$ admits 
    infinitely many submodels isomorphic to~$M$.
    \item Any isomorphism between two finite submodels 
    of $\mathbf{F}$ can be extended to an automorphism 
    of $\mathbf{F}$.
\end{itemize}
\end{prp}

\begin{proof}
The construction of $\mathbf{F}$ is made by induction. 
Let $F_0=\emptyset$. $F_0$ is trivially a model of $\T m$ and, 
for every element $M$ of $C$ with $|M|\leq0$, $F_0$ contains a 
submodel isomorphic to $M$. Let $F_n$ be a model of $\T m$ 
with a finite initial segment of $\omega$ as domain and 
such that for every $M\in C$ with $|M|\leq n$, $F_n$ 
contains a submodel isomorphic to $M$. Let
\begin{itemize}
    \item $\{A_i:i\leq p\}$ be an enumeration of $[F_n]^{\leq n}$,
    \item $\{R_k: k\leq q\}$ be an enumeration of all the 
    $M\in C$ such that $1\leq|M|\leq n+1$,
    \item $\{j_l:l\leq u\}$ be an enumeration of all the 
    embeddings $j_l:F_n|_{A_i}\xhookrightarrow{} R_k$, 
    where $i\leq p$, $k\leq q$ and $|R_k|=|A_i|+1$.
\end{itemize}
For each $l\leq u$, let $a_l\in\omega$ be the least natural 
number such that $a_l\notin F_n\cup\{a_{l'}:l'<l\}$. 
The idea is to add $a_l$ to $F_n$, extending $F_n|_{A_i}$ to a 
model $F_n|_{A_i}\cup\{a_l\}$ isomorphic to $R_k$, 
where $j_l:F_n|_{A_i}\xhookrightarrow{} R_k$. 
Define $F_{n+1}:=F_n\cup\{a_l:l\leq u\}$. In \cite{cst}, 
$F_{n+1}$ is made into a model of $\T m$ in a non-controlled way, 
while here we impose the following: 
Let $\{x_1,\dots,x_{n}\}$ be a subset of $F_{n+1}$ from which 
we have not already chosen an element. Suppose that $s>t$ 
implies $x_s>x_t$ (recall that $F_{n+1}$ is a subset of $\omega$). 
Then we simply impose $\Sel{}(\{x_1,\dots,x_n\})=x_n$. 
The desired model is finally given by 
$\mathbf{F}=\bigcup_{n\in\omega}F_n$.

We conclude by showing that every isomorphism between finite 
submodels can be extended to an automorphism of $\mathbf{F}$. 
Let $i_0:M_1\to M_2$ be an isomorphism of $\T m$-models. 
Let $a_1$ be the least natural number in 
$\omega\setminus(M_1\cup M_2)$. 
Then $M_1\cup M_2\cup\{a_1\}$ is contained in some $F_n$ 
and by construction we can find some $a'_1\in\omega$ such 
that $\mathbf{F}|_{M_1\cup\{a_1\}}$ is isomorphic to 
$\mathbf{F}|_{M_2\cup\{a'_1\}}$. Extend $i_0$ to 
$i_1:M_1\cup\{a_1\}\to M_2\cup\{a'_1\}$ by imposing 
$i_1(a_1)=a'_1$. Let $a_2$ be the least integer in 
$\omega\setminus(M_1\cup M_2\cup\{a_1,a'_1\})$ and repeat 
the process. The desired automorphism of $\mathbf{F}$ 
is $i=\bigcup_{n\in\omega}i_n$. 
\end{proof}

\begin{rmk}
Let us fix some notations and terminology. The elements 
of the model $\mathbf{F}$ above constructed will be the 
atoms of our permutation model. Each element $a$ corresponds 
to a unique embedding $j$. We shall call the domain of $j$ 
the $\emph{ground}$ of $a$. Moreover, given two atoms $a$ and $b$, 
we say that $a<b$ in case $a<_\omega b$ according to 
the natural ordering. Notice that this well ordering 
of the atoms will not exist in the permutation model.
\end{rmk}

Let $A$ be the domain of the model $\mathbf{F}$ of the 
theory $\T m$. To build the permutation model $\MOD m$, 
consider the normal ideal given by all the finite subsets 
of $A$ and the group of permutations $G$ defined by
$$\pi\in G\iff \forall\,X\in\fin(\omega),
\pi(\Sel{}(X))=\Sel{}(\pi X).$$

\begin{thm}\label{thm:7.2}
For every positive integer $m$, $\MOD m$ is a model for $\RC m$.
\end{thm}

\begin{proof}
Let $X$ be an infinite set with support $S'$. If $X$ is well 
ordered, the conclusion is trivial, so let $x\in X$ be an element 
not supported by $S'$ and let $S$ be a support of $x$, 
with $S'\subseteq S$. Let $a\in S\setminus S'$. 
If $\fix_G(S\setminus\{a\})\subseteq\sym_G(x)$ then 
$S\setminus\{a\}$ is a support of $x$, so by iterating the 
process finitely many times we can assume that there exists 
a permutation $\tau\in\fix_G(S\setminus\{a\})$ such that 
$\tau(x)\neq x$. Our conclusion will follow by showing that 
there is a bijection between an infinite set of atoms and 
a subset of $X$, namely between 
$\{\pi(a):\pi\in\fix_G(S\setminus\{a\})\}$ and 
$\{\pi(x):\pi\in\fix_G
(S\setminus\{a\})\}$. Suppose towards a contradiction that 
there are two permutations
 $\sigma,\sigma'\in\fix_G(S\setminus\{a\})$ such that 
 $\sigma(x)=\sigma'(x)$ and $\sigma(a)\neq\sigma'(a)$. 
 Then, by direct computation, the permutation $\sigma^{-1}\sigma'$ 
 is such that $\sigma^{-1}\sigma'(a)\neq a$ and 
 $\sigma^{-1}\sigma'(x)=x$. Let $b=\sigma^{-1}\sigma'(a)$. 
 Then $\{b\}\cup (S\setminus\{a\})$ is a support of $x$. 
 By construction of the model, the set 
 $\{\pi(a):\pi\in\fix_G(\{b\}\cup (S\setminus\{a\}))\}$ 
 is infinite, from which we deduce that also the set 
$$
L=\{e\in A: \exists\,\pi\in\fix_G(\{b\}\cup(S\setminus\{a\}))
\textrm{ such that }\pi(x)=x,\,\pi(a)=e\}
$$
is infinite. Now, by assumption
there is a permutation $\tau\in\fix_G(S\setminus\{a\})$ 
such that $\tau(x)\neq x$. Let $y_0:=\tau(x)$, with 
$c=\tau(a)$ and $d=\sigma^{-1}\sigma'(c)$. Then the same 
argument shows that also 
$$
R=\{e\in A: \exists\,\pi\in\fix_G(\{d\}\cup(S\setminus\{a\}))
\textrm{ such that }\pi(x)=y_0,\,\pi(c)=e\}
$$ 
must be infinite.\newline

\noindent First note that in $L$ there are infinitely many 
elements with ground $\{b\}\cup(S\setminus \{a\})$, which 
is a support of $L$. This is because 
$\{b\}\cup(S\setminus \{a\})\cup \{a\}\subseteq \textbf{F}$ 
is a finite model of $\T m$ and in the construction of 
our permutation model we add infinitely many elements 
$l\in\omega$ such that $\{b\}\cup(S\setminus \{a\})
\cup \{l\}$ and $\{b\}\cup(S\setminus \{a\})\cup \{a\}$ 
are isomorphic via an isomorphism $\delta$ with 
$\delta\vert_{\{b\}\cup(S\setminus \{a\})}=
\operatorname{id}\vert_{\{b\}\cup(S\setminus \{a\})}$ 
and $\delta(l)=a$. We can extend $\delta$ to an 
automorphism $\delta\in\fix_G(\{b\}\cup(S\setminus \{a\}))$. 
By definition of $L$, since $\{b\}\cup(S\setminus \{a\})$ 
is a support of $x$, we have that $l\in L$.

Similarly, in $R$ there are infinitely many elements 
with ground $\{d\}\cup(S\setminus \{a\})$, which is a support 
of $R$. This is because $\{d\}\cup(S\setminus \{a\})\cup 
\{c\}\subseteq \textbf{F}$ is a finite model of $\T m$ 
and in the construction of our permutation model we add 
infinitely many elements $l\in\omega$ such that 
$\{d\}\cup(S\setminus \{a\})\cup \{l\}$ and 
$\{d\}\cup(S\setminus \{a\})\cup \{c\}$ are isomorphic 
via an isomorphism $\delta$ with $\delta\vert_{\{d\}
\cup(S\setminus \{a\})}=\operatorname{id}\vert_{\{d\}
\cup(S\setminus \{a\})}$ and $\delta(l)=c$. We can 
extend $\delta$ to an automorphism $\delta\in
\fix_G(\{d\}\cup(S\setminus \{a\}))$. By definition 
of $R$, since $\{c\}\cup(S\setminus \{a\})$ is a 
support of $y_0$, we have that $l\in R$.
\newline

\noindent Let $r\in R$ with ground $\{d\}\cup(S\setminus 
\{a\})$ and $p,l\in L$ with ground $\{b\}\cup(S\setminus
\{a\})$ such that $r\geq p$, $l\geq p$ and 
$\min(\{p,q,r\})>\max(\{d,b\}\cup(S\setminus \{a\}))$. 
We want to show that every map
$$
\gamma:(S\setminus \{a\})\cup \{p\}\cup\{l\}\to 
(S\setminus \{a\})\cup \{p\}\cup \{r\}
$$
with $\gamma\vert_{(S\setminus \{a\})\cup \{p\}}=
\operatorname{id}_{(S\setminus \{a\})\cup \{p\}}$ and 
$\gamma(l)=r$ is an isomorphism of $\T m$-models. Let 
$X\subseteq (S\setminus \{a\})\cup \{p\}\cup\{l\}$. 
If $\{p,l\}\cap X=\emptyset$ we have that $\gamma(\Sel{}(X))=
\Sel{}(\gamma(X))$. If $l\in X$ and $p\notin X$ let 
$\pi_l,\pi_r\in\fix_G(S\setminus\{a\})$ with $\pi_l(a)=l$ and 
$\pi_r(a)=r$. Then $\pi_r\circ\pi_l^{-1}\vert_X=\gamma\vert_X$. 
So since $\pi_r\circ\pi_l^{-1}\in G$ we have 
$\gamma(\Sel{}(X))=\Sel{}(\gamma(X))$. In the last case, 
when $\{p,l\}\subseteq X$, the function selects the greatest 
element with respect to $<_\omega$, given the particular care 
we took in the construction of the selection function on the 
set of atoms. So we can extend $\gamma$ to a function 
$\tau^{\prime}\in\fix_G(\{p\}\cup (S\setminus \{a\}))$ 
with $\tau^{\prime}(l)=r$.\newline

\noindent Let $\pi_r\in\fix_G(S\setminus\{a\})$ such that 
$\pi_r(a)=r$ and $\pi_r(x)=y_0$. 
Let $\pi_l\in\fix_G(S\setminus\{a\})$ with $\pi_l(a)=l$ 
and $\pi_l(x)=x$. Then we have that 
$\pi_r^{-1}\circ\tau^{\prime}\circ\pi_l(a)=a$ which 
implies that $\pi_r^{-1}\circ\tau^{\prime}\circ\pi_l(x)=x$ 
because the function fixes $S$. So 
\begin{equation}
\label{eq:1}
\tau^{\prime}(x)=\tau^{\prime}\circ\pi_{l}(x)=\pi_r(x)=y_0.
\end{equation}
Now let $\pi_p\in \fix_G(S\setminus\{a\})$ with $\pi_p(a)=p$ 
and $\pi_p(x)=x$. Since $S$ is a support of $x$, 
$\pi_p(S)=\{p\}\cup (S\setminus\{a\})$ is also a support 
of $\pi_p(x)=x$. Therefore,
$$
\tau^{\prime}(x)=x.
$$
This is a contradiction to (\ref{eq:1}). So we showed that 
for all $\sigma,\sigma'\in\fix_G(S\setminus\{a\})$, 
$\sigma(x)=\sigma'(x)$ implies $\sigma(a)=\sigma'(a)$, 
from which we get the desired bijection.
\end{proof}

\section{For which $\boldsymbol n$ is $\boldsymbol{\MOD m}$ 
a model for $\boldsymbol{\RC n}$\,?}

The following result shows that for positive integers $m,n$ 
which satisfy a certain condition,
the implication $\RCmn mn$ is not provable in $\ZF$. 
Assuming the ternary Goldbach conjecture, it will 
turn out that all positive integers $m,n$ 
satisfy this condition, 
except when $m=n$, or when $m=2$ and~$n=4$.

\begin{prp}\label{prp:negative_general}
    Let $m,n\in\omega$. Assume that, for some $k\in\omega$, 
    $n$ can be written as a sum $n=\sum_{i\in k}n_i$ 
    with each $n_i\neq 1$, such that for all possible 
    finite sequences $(m_i)_{i\in k}$ satisfying 
\begin{itemize}
        \item $0\leq m_i\leq n_i$
        \item $m_i = 0$ or $\gcd(n_i,m_i)>1$
\end{itemize}
    we have that $m\neq\sum_{i\in k} m_i$. 
    Then the implication $\RCmn mn$ is not provable in $\ZF$.
\end{prp}

\begin{proof}
    We show that in $\MOD m$, $\RC n$ fails. Assume 
    towards a contradiction that $\RC n$ holds in 
    $\MOD m$ and let $S$ be a support of a selection 
    function $f$ on the $n$-element subsets of an 
    infinite subset $X$ of the set of atoms $A$.
    
    By the construction in Proposition\;\ref{prp:2.1}, given any model 
    $N$ of $\T m$ extending $S$, we can find a submodel 
    of $X\cup S$ isomoprhic to $N$. 
    
    Our conclusion can hence follow from finding a model 
    $M$ of $\T m$ 
    which extends $S$ with $|M\setminus S|=n$ and such that 
    $M$ admits an auotmorphism $\sigma$ which fixes pointwise 
    $S$ and which does not have any other fixed point, 
    since then $\sigma(f(M\setminus S))\neq f(M\setminus S)$ 
    but $\sigma(M\setminus S) = M\setminus S$. We start 
    with the following claim:

    \begin{claim}
        Given a cyclic permutation $\pi$ on some set $P$ 
        of cardinality $|P|=q$, if a non-trivial power 
        $\pi^r$ of $\pi$ fixes a proper subset $P'$ of 
        $P$, then $\gcd(|P'|,|P|)>1$.
    \end{claim}
    
    To prove the claim, notice that $\pi^r$ is a disjoint 
    union of cycles of the same length $l=\frac{q}{\gcd(q,r)}$. 
    Consider the subgroup of $\langle\pi\rangle$ given 
    by $\langle\pi^r\rangle$. Then $P'$ is a disjoint 
    union of orbits of the form $\textrm{Orb}_{<\pi^r>}(e)$ 
    with $e\in P'$, all of them with the same cardinality $s$, 
    with $s$ being a divisor of $l=\frac{q}{\gcd(q,r)}$ 
    and hence of $q$, from which we deduce the claim.
     
We now want to show that we can find a model $M$ of $\T m$, 
which extends $S$ with $|M\setminus S|=n$ and such that it admits 
an automorphism $\sigma$ which fixes pointwise $S$ and acts on 
$M\setminus S$ as a disjoint union of $k$ cycles, each of length 
$n_i$ for $i\in k$. This can be done as follows. Pick an $m$-element 
subset $P$ of $M$ for which $\Sel{}(P)$ has not been defined yet. I
f $P\cap S\neq\emptyset$ then let $\Sel{}(P)$ be any element 
in $P\cap S$. Otherwise, by our the assumptions, 
there is a cycle $C_j$ of length $n_j$ for some $j\in k$ 
such that $\gcd(|P\cap C_j|, |C_j|) = 1$. Define $\Sel{}(P)$ 
as an arbitrarily fixed element of $P\cap C_j$ and, for all 
permutations $\pi$ in the group generated by $\sigma$, define 
$\Sel{}(\pi(P))=\pi(\Sel{}(P))$. We need to argue that this 
is indeed well defined, i.e. that for two permutations 
$\pi,\pi'\in\langle\sigma\rangle$ we have that $\pi(P)=\pi'(P)$ 
implies $\pi(\Sel{}(P))=\pi'(\Sel{}(P))$. 
Problems can arise only when $P\cap S=\emptyset$, in which 
case we notice that $\pi(P)=\pi'(P)$ implies $\pi(P\cap C_j)=
\pi'(P \cap C_j)$, which in turn by the claim implies that 
$\pi^{-1}\circ\pi'$ fixes $P\cap C_j$ pointwise, from 
which we deduce $\pi(\Sel{}(P))=\pi'(\Sel{}(P))$.
\end{proof}

Proposition\;\ref{prp:negative_general} allows us to 
immediately deduce the following results.

\begin{cor}\label{greater}
    If $m>n$, then $\RC m$ does not imply $\RC n$.
\end{cor}

\begin{proof}
    It is enough to write $n$ as $n= \sum_{i\in 1}n_i$ 
    with $n_0=n$ and apply 
    Proposition\;\ref{prp:negative_general}.
\end{proof}

\begin{cor}\label{primes_inclusion}
    If there is a prime $p$ for which $p\mid n$ but $p\nmid m$, 
    then $\RC m$ does not imply $\RC n$.
\end{cor}

\begin{proof}
    Given the assumption, writing $n$ as $n=\sum_{i\in 
    \frac{n}{p}}n_i$, where each $n_i=p$, allows us to 
    directly apply Proposition\;\ref{prp:negative_general}.
\end{proof}

Moreover, we can show the following:

\begin{thm}\label{thm:main}
For any positive integers\/ $m$ and\/ $n$, 
the implication\/ $\RCmn mn$ is provable in\/ $\ZF$ only in
the case when\/ $m=n$ or when~$m=2$ and~$n=4$.
\end{thm}

The proof of Theorem\;\ref{thm:main} is given in the 
following results, where in the proofs we  
use two well-known number-theoretical 
results: The first one is Bertrand's postulate, 
which asserts that for every positive
integer $m\ge 2$ there is a prime $p$ with $m<p< 2m$, 
and the second one is ternary Goldbach 
conjecture (assumed to be proven by 
Helfgott~\cite{goldbach}), which asserts that every 
odd integer $n>5$ is the sum of three primes.   

\begin{prp}\label{prp:prime}
    If\/ $m$ is prime and\/ $n\neq m$, then the 
    implication\/ $\RCmn mn$ is not provable in\/ $\ZF$
\end{prp}

\begin{proof}
    Given Corollary\;\ref{primes_inclusion}, we can assume that 
    $n = m^k$ for some natural number $k>1$. Let $p$ be a prime 
    such that $m<p<2m$, whose existence is guaranteed by 
    Bertrand's postulate. Then clearly $m\nmid n-p$, from 
    which we get that the decomposition $n = p + (n-p)$ 
    satisfies the assumptions of 
    Proposition\;\ref{prp:negative_general}.
\end{proof}

\begin{prp}\label{prp:odd}
    If\/ $n$ is odd and\/ $m\neq n$, then the implication\/
    $\RCmn mn$ is not provable in\/ $\ZF$.
\end{prp}

\begin{proof}
    By the ternary Goldbach conjecture, let us write 
    $n$ as sum of three primes $n=p_1+p_2+p_3$. 
    Given Proposition\;\ref{prp:prime}, we can assume that $m=p_1+p_2$, 
    since otherwise the decomposition $n=p_1+p_2+p_3$ would 
    satisfy the assumptions of Proposition\;\ref{prp:negative_general}.

    We first deal with the case in which $p_1=p_2=p_3$ holds, 
    for which we rename $p=p_1$. By hand we can exclude the 
    case $p=2$, and now we want to show that the decomposition 
    $n=n_0+n_1= (3p - 2) + 2$ can be used in 
    Proposition\;\ref{prp:negative_general}. Notice that 
    $\gcd(3p-2,2p-2)\in\{1,p\}$, from which we deduce that 
    necessarily if $m=m_0+m_1$ as in 
    Proposition\;\ref{prp:negative_general}, then $m_1 = 0$. 
    To conclude this first case, it suffices to notice that, 
    since $p$ is a prime grater than $2$, $\gcd(3p-2, 2p)$ 
    necessarily equals $1$.

We can now assume that it is not true that $p_1=p_2=p_3$. Since $n$ is odd, 
$p_1+p_2\nmid p_3$. If $p_3\nmid p_1+p_2$, then the decomposition $n=n$ 
actually satisfies Proposition\;\ref{prp:negative_general}. So, without
loss of generality let us assume that $p_3<p_1$ and $p_3\mid p_1+p_2$. 
By $p_3\mid p_1+p_2$ 
we deduce that $p_2\neq p_3$, and we now consider the decomposition 
$n=n_0+n_1=(p_2+p_3) + p_1$. We can't have $m_1=p_1$ since 
$\gcd(p_2,p_2+p_3)=1$. On the other hand, we can't even have 
$m_1=0$ since $p_1+p_2 > p_2+p_3$, which proves that the 
assumptions of Proposition\;\ref{prp:negative_general} are satisfied.
\end{proof}

\begin{prp}\label{prp:p2}
    Let\/ $m>2$ be an even natural number and\/ $k\in\omega$ 
    such that\/ $2^k+1$ is prime. If\/ $n=m+2^k$, then the implication\/ 
    $\RCmn mn$ is not provable in\/ $\ZF$.
\end{prp}

\begin{proof}
We consider the decomposition $n=n_0+n_1=(m-1)+(2^k+1)$. 
It directly follows from the assumptions of the proposition that in 
order to write $m=m_0+m_1$ as in Proposition\;\ref{prp:negative_general}, 
since $n_0<m$, necessarily $m_1=2^k+1$, from which we deduce 
$m_0=m-2^k-1$. This immediately gives a contradiction in the 
case $2^k+1>m$, so let us assume $2^k+1<m$. We get again a 
contradiction by the fact that 
$\gcd(m_0,n_0)=\gcd(m-2^k-1,m-1)=\gcd(2^k,m-1)=1$, where we 
used that $m$ is even. We can hence conclude that 
Proposition\;\ref{prp:negative_general} can be applied.
\end{proof}

\begin{prp}\label{prp:even1}
Let\/ $m$ and\/ $n$ be even natural numbers such that there 
is an odd prime\/ $p$ with\/ $m<p<n$ and\/ $n>p+1$. 
Then the implication\/ $\RCmn mn$ is not provable in\/ $\ZF$.
\end{prp}

\begin{proof}
If $n=p+3$ or $n=p+5$ the decomposition $n= p +(n-p)$ directly allows us to apply Proposition\;\ref{prp:negative_general}. Otherwise, by the ternary Goldbach conjecture, 
write $n-p$ as sum of three primes $n-p=p_1+p_2+p_3$. 
Consider now the decomposition $n=\sum_{i\in4}n_i=p+p_1+p_2+p_3$. 
In order to write $m=\sum_{i\in4}m_i$, necessarily $m_0=0$. 
If $n-p<m$ we can already conclude that 
Proposition\;\ref{prp:negative_general} can be applied. 
Otherwise, we find ourselves in the assumptions of 
Proposition\;\ref{prp:odd}, which again allows us to conclude that 
$\RC m$ does not imply $\RC n$.
\end{proof}

The following result deals with all the remaining cases and 
completes the proof of Theorem\;\ref{thm:main}.

\begin{prp}\label{prp:even2}
Let\/ $m$ and\/ $n$ be even natural numbers with\/ $3\leq\frac{n}{2}\leq m<n$ 
such that if there is a prime\/ $p$ with\/ $m<p<n$, then $p=n-1$. 
Then the implication\/ $\RCmn mn$ is not provable in\/ $\ZF$.
\end{prp}

\begin{proof}
By Bertrand's postulate, let $p$ be a prime with 
$\frac{n}{2}<p<n$. This implies by the assumption $\frac{n}{2}<p<m$ or $p=n-1$. If we are in the latter case, apply again Bertrand's postulate 
to find a further prime $\frac{n}{2}-1<p'<n-2$
(notice that by our assumption we have
$2\leq \frac{n}{2}-1$). Since $m$ is not prime we necessarily have $p'\neq m$, which together with the present assumptions makes us able to assume without loss of generality that $\frac{n}{2}<p<m$.
By Proposition\;\ref{prp:p2} we can assume $n-m>4$, which implies $n-p>5$.
Since by the ternary Goldbach conjecture we can write  we 
can write $n=p+p_1+p_2+p_3$, notice that by the fact that 
$n$ and $m$ are even, we can assume that $m-p$ equals some 
odd prime $p'$, since otherwise we could already apply 
Proposition\;\ref{prp:negative_general} with the decomposition 
$n=p+p_1+p_2+p_3$ or Proposition\;\ref{prp:p2} with the same pair $(m,n)$. 
Now, either writing $n=p+(n-p)$ satisfies the assumptions of 
Proposition\;\ref{prp:negative_general}, or $n-p$ is a multiple 
of $p'$. We distinguish two cases, namely when $n-p$ is a 
power of $p'$ and when it is not. In the second case, let 
$p''$ be a prime distinct from $p'$ such that $p''\mid n-p$. 
Writing $n$ as the sum 
$n=n_0+\sum_{i\in\frac{n-p} {p''}}n_i=p+\sum_{i\in\frac{n-m}{p''}}p''$ 
satisfies the hypothesis of Proposition\;\ref{prp:negative_general}, 
as $n-p<m$ and hence if $m=m_0+\sum_{i\in\frac{n-m}{p'}}m_i$ 
then $m_0=p$. For the last case, without loss of generality assume
that $p_1+p_2+p_3=p_1^k$ 
for some natural number $k>1$. Notice that $p_3\neq p_1$, 
since otherwise we would have $p_2 = p_1^k-2p_1=p_1(p_1^{k-1}-2)$, 
which is a contradiction, and similarly we obtain $p_2\neq p_1$. 
We finally assume wlog that $p_2>p_1$, which allows us to 
conclude that the sum $n=p+p_2+(p_1+p_3)$ satisfies the 
assumptions of Proposition\;\ref{prp:negative_general}, 
concluding the proof.
\end{proof}

For the sake of completeness, we summarise the proof of our main theorem:

\begin{proof}[Proof of Theorem\;\ref{thm:main}]
Let $m$ and $n$ be two distinct positive integers.
$$\ZF\proves \RCmn mn
\quad\overset{\text{Cor.\,3.2}}{\Longrightarrow}\quad m<n
\quad\overset{\text{Prp.\,3.6}}{\Longrightarrow}\quad \text{$n$ is even}
\quad\overset{\text{Cor.\,3.3}}{\Longrightarrow}\quad \text{$m$ is even}$$
Now, if $m$ and $n$ are both even, we have the following two cases:
\begin{eqnarray*}
m<\frac n2 
& \overset{\text{Prp.\,3.8}}{\Longrightarrow} 
& \ZF\not\proves \RCmn mn\\[1ex]
m\geq\frac n2\geq 3 
& \underset{\text{Prp.\,3.8}}{\overset{\text{Prp.\,3.9}}{\Longrightarrow}}
& \ZF\not\proves \RCmn mn
\end{eqnarray*}
Thus, by Fact\;\ref{fct:RC24}, the implication $\RCmn mn$ is 
provable in $\ZF$ if and only if $m=n$ or $m=2$ and $n=4$.
%
%
\end{proof}

\end{document}